\documentclass{article}
\usepackage[utf8]{inputenc}
\usepackage{amsthm, amssymb}
\usepackage{amsmath}

\usepackage[normalem]{ulem} 

\title{Cardinality of definable families of sets in o-minimal structures}
\author{Pablo And\'ujar Guerrero \\
\small \emph{Email address:} pa377@cantab.net}
\date{}

\newtheorem{definition}{Definition}
\newtheorem{lemma}{Lemma}
\newtheorem{fact}{Fact}
\newtheorem{proposition}{Proposition}

\newtheorem{claim}{Claim}
\newtheorem{corollary}{Corollary}

\newcommand{\Dd}{\mathcal{D}}
\newcommand{\Cc}{\mathcal{C}}
\newcommand{\df}{\mathrm{def}}
\newcommand{\cl}{cl}
\newcommand{\infM}{M \cup \{-\infty, +\infty\}}
\newcommand{\Ff}{\mathcal{F}}
\newcommand{\Gg}{\mathcal{G}}
\newcommand{\Xx}{\mathcal{X}}
\newcommand{\Def}{\mathrm{Def}}

\newenvironment{claimproof}[1][\proofname]
               {
                 \proof[#1]
                 
               }
               {
                 \endproof
               }

\begin{document}
\maketitle

\noindent
{\small \emph{2020 Mathematics Subject Classification:} 03C64 (Primary), 03C45 (Secondary).} 

\begin{abstract}
We prove that any definable family of subsets of a definable infinite set $A$ in an o-minimal structure has cardinality at most $|A|$. We derive some consequences in terms of counting definable types and existence of definable topological spaces.      
\end{abstract}

\section{Introduction} 

Let $T$ be a theory. A (partitioned) formula $\varphi(x;y)$ is stable\footnote{For a more extensive characterization of stability and its connection to dependence see~\cite{NIPguide}: Proposition 2.55, Definition 2.56 and the discussion in Section 2.3.4} (in $T$) if there exists some $k<\omega$ such that, for any model $M\models T$ and set $A\subseteq M^{|x|}$, it holds that 
\begin{equation}\label{eqn:intr}
|\{\varphi(A;b) : b\in M^{|y|} \}| \leq |A|^{k}.
\end{equation}
And in the specific case where $A$ is infinite  
\begin{equation}\label{eqn:intr0}
|\{\varphi(A;b) : b\in M^{|y|} \}| \leq |A|.
\end{equation}
In the terminology of types from~\cite{NIPguide}, the set of (complete) $\varphi^{\text{opp}}$-types over $A$ is bounded in size by $|A|^k$. A theory is stable if every formula is stable. The theory of any infinite linear order is unstable (not stable). In particular for $M=(\mathbb{R},<)$ consider the formula $x<y$ and set $A=\mathbb{Q}$.

NIP theories are a generalization of stable theories that are characterized by the property that, for every formula $\varphi(x;y)$, there exists some $k<\omega$ such that inequality \eqref{eqn:intr} holds for every \textbf{finite} set $A$.
Meanwhile, whenever $A$ is infinite, inequality \eqref{eqn:intr0} fails in general, even with the assumption that $A$ is definable in $M$, as witnessed by the structure $(\mathbb{R}, <, \mathbb{Q})$, which is NIP by~\cite[Corollary 3.2]{gun-hie-dependent}.


In this note we show (Proposition~\ref{prop:main}) that in o-minimal theories, a subclass of NIP theories, inequality~\eqref{eqn:intr0} holds whenever $A$ is an infinite set \textbf{definable} in $M$. In the terminology of types, for every formula $\varphi(x;y)$ the family of realized $\varphi$-types over some infinite definable set $A$ is bounded in size by $|A|$.  

As a corollary of our result, we establish bounds on the size of various spaces of definable types over $A$. We also bound the topological weight of any o-minimal definable topological space, showing in particular that any countable definable topological space is second-countable. 

It is worth noting that our result is easy to prove in o-minimal expansions of ordered fields, using that any two intervals are in bijection, and so by o-minimality any two infinite definable sets have the same cardinality.

\section{Cardinality of definable families of sets}

\subsection{Conventions}\label{sec:dfns}

Throughout we work in a structure denoted $M$. By ``definable" we mean ``definable in $M$ possibly with parameters".  Any formula we consider is in the language of $M$.
We use $x$, $y$, $z$ and $a$, $b$, $c$ to denote tuples of variables and tuples of parameters respectively. We use $s$ and $t$ exclusively for unary parameters.

Let $\varphi(x)$ be a formula, possibly with parameters from $M$. For any set $A\subseteq M^{|x|}$ let $\varphi(A)=\{ a\in A : M\models \varphi(a)\}$. We write $\varphi(M)$ to mean $\varphi(M^{|x|})$.

For background o-minimality we direct the reader to~\cite{dries98}. We will use in particular o-minimal cells, the existence of uniform cell decompositions~\cite[Chapter 3, Proposition 3.5]{dries98} and the Fiber Lemma for o-minimal dimension~\cite[Chapter 4, Proposition 1.5]{dries98}. 
We also adopt from~\cite{dries98} the convention of the unary space $M^0$, where any parameter $a\in M^n$ is implicitly understood to have its ``zero coordinate" in $M^0$, and moreover $a$
can be seen as a function $M^0\rightarrow M^n$ with image $\{a\}$.

For $M$ o-minimal and $A\subseteq M$ a set we denote by $cl(A)$ the closure of $A$ in the o-minimal order topology. We also extend the linear order in $M$ in the natural way to $\infM$. For a function $f$ we denote its domain by $dom(f)$ and its graph by $graph(f)$. We adopt the usual convention surrounding o-minimal cells of saying that a partial function $M^{n}\rightarrow \infM$ is definable if either it maps into $M$ and is definable in the usual sense or otherwise it is constant and its domain is definable. For a set $B\subseteq M^{n}$, we say that a family $\{f_b : b\in B\}$ of definable partial functions $M^{n}\rightarrow \infM$ is definable if the sets $B(+\infty)=\{ b\in B : f_b \equiv +\infty\}$ and $B(-\infty)=\{ b\in B : f_b \equiv -\infty\}$ are both definable and moreover the families $\{ dom(f_b) : b\in B(+\infty)\}$, $\{ dom(f_b) : b\in B(-\infty)\}$ and $\{f_b : b\in B \setminus (B(+\infty)\cup B(-\infty))\}$ are (uniformly) definable in the usual sense.


In general types are complete and consistent. 
A type $p(x)$ over a set $A\subseteq M$ is definable if, for every formula $\varphi(x;y)$, there is another formula $\psi(y;z)$ and some $c\in M^{|z|}$ such that 
\[
\{b \in M^{|y|} : \psi(x;b) \in p(x) \} = \psi(M;c).
\]
Notably this differs from other literature definitions that require that $c\in A^{|z|}$. We denote by $S^{\df}(A)$ the space of types over $A$ that are definable. We denote by $S_x^{\df}(A)$ the subspace of $S^{\df}(A)$ of types $p(x)$ with object variable $x$. 

\subsection{Main result and corollaries}

We begin with a simple lemma. 

\begin{lemma}\label{lem:size}
Let $\Xx$ be a family of sets and $\{ \Xx_i : i \leq k \}$ be a finite collection of families of sets such that, for every $A \in \Xx$, there exists $I\subseteq \{1,\ldots,k\}$ and $\{A(i) \in \Xx_i : i \in I\}$ such that $A= \cup_{i \in I} A(i)$. Then 
\[
|\Xx| \leq \max_{i \leq k} |\Xx_i|+\aleph_0. 
\]
\end{lemma}
\begin{proof}
We may assume that, for every $i\leq k$, the family $\Xx_i$ contains the empty set. Hence, for every $A\in \Xx$, the set $I$ in the lemma can always be taken to be $\{1,\ldots, k\}$. 

Observe that any function
\begin{align*}
\Xx &\rightarrow \Xx_1 \times \cdots \times \Xx_k \\
 A &\mapsto (A(1), \ldots, A(k)),
\end{align*}
where $A= \cup_{i\leq k} A(i)$, is injective. Hence 
\[
|\Xx| \leq \prod_{i\leq k}  |\Xx_i| \leq \max_{i\leq k} |\Xx_i| + \aleph_0. 
\]
\end{proof}

We now present our main result. 

\begin{proposition}\label{prop:main}
Let $M$ be an o-minimal structure. For any formula $\varphi(x;y)$ and infinite definable set $A \subseteq M^{|x|}$ it holds that 
\[
|\{\varphi(A;b) : b\in M^{|y|}\}| \leq |A|. 
\]
\end{proposition}
\begin{proof}
Fix $A \subseteq M^{|x|}$ and $\varphi(x;y)$ as in the proposition. We proceed by induction on $|x|+|y|$. 
We begin the proof by reducing it using o-minimal cell decomposition to the case where, for every $b\in M^{|y|}$, the formula $\varphi(A;b)$ defines the graph of some partial function $M^{|x|-1}\rightarrow M$. 


Throughout let $\pi:M^{|x|}\rightarrow M^{|x|-1}$ be the projection to the first $|x|-1$ coordinates. For every $c\in M^{|x|-1}$ let $A_c$ denote the fiber $\{ t : (c,t) \in A\}$. Let $\overline{A}:=\{ (c,t) : c\in \pi(A),\, t\in \cl(A_c)\}$. Observe that this set is definable.
By o-minimality there exists $n < \omega$ such that, for every $c\in M^{|x|-1}$, it holds that $\cl(A_c) \leq |A_c|+n$, and so $|\overline{A}|=|A|$. 

By o-minimal uniform cell decomposition~\cite[Chapter 3, Proposition 3.5]{dries98} there exist finitely many formulas $\psi_1(x;y), \ldots, \psi_k(x;y)$ such that, for any $b\in M^{|y|}$, the family $\{ \psi_1(M;b), \ldots, \psi_k(M;b)\}$ is an o-minimal cell decomposition of $M^{|x|}$ compatible with $\{\varphi(A;b), A\}$. By Lemma~\ref{lem:size}, the proposition follows from showing that 
\begin{equation}\label{eqn:1}
|\{\psi_i(A;b) : b\in M^{|y|}\}| \leq  |\overline{A}| = |A|
\end{equation}
for every $i\leq k$. 

Fix $i\leq k$. For any $b \in M^{|y|}$, let $f_{b}$ and $g_b$ be the definable partial functions $M^{|x|-1}\rightarrow \infM$ such that either $\psi_i(M;b)=(f_b,g_b)$ or $\psi_i(M;b)=graph(f_b)=graph(g_b)$. Let $B = \{ b\in  M^{|y|} : \psi_i(A,b)\neq \emptyset\}$. Note that the families $\{ f_b : b\in B\}$ and $\{ g_b : b\in B\}$ are definable, and moreover \eqref{eqn:1} follows from showing that each one has cardinality at most $|\overline{A}|$. 
Now since, for every $b\in M^{|y|}$, the set $\psi_i(M,b)$ is either disjoint or contained in $A$, we have that, for any $b\in B$, the graphs of $f_b$ and $g_b$ are contained in $\overline{A} \cup (\pi(A) \times \{-\infty, +\infty\})$.

Consider the definable families of formulas $\Ff_{\infty}=\{ f_b : f_b \equiv -\infty,\, dom(f_b)\subseteq \pi(A),\, b\in M^{|y|}\}$ and $\Gg_{\infty}=\{ g_b : g_b \equiv +\infty,\, dom(g_b)\subseteq \pi(A),\, b\in M^{|y|}\}$. By induction hypothesis (the case $|x|=1$ being trivial) applied to the definable families $\{dom(f_b) : f_b\in \Ff_{\infty}\}$ and $\{dom(g_b) : g_b\in \Gg_{\infty}\}$, we have that $\Ff_{\infty}$ and $\Gg_{\infty}$ are both bounded in size by $|A|=|\overline{A}|$.
Hence to prove \eqref{eqn:1} it suffices to bound the size of $\Ff_M=\{ f_b : graph(f_b) \subseteq \overline{A},\, b\in M^{|y|}\}$ and $\Gg_M=\{ g_b : graph(g_b)\subseteq \overline{A},\, b\in M^{|y|}\}$ by $|\overline{A}|$. 

We conclude that, by passing if necessary and without loss of generality from $A$ and $\{\varphi(A;b) : b\in M^{|y|}\}$ to $\overline{A}$ and $\Ff_M$ (or $\Gg_M$) respectively, we may assume that the formula $\varphi(x;y)$ is such that, for any $b\in M^{|y|}$, the set $\varphi(A;b)$ defines the graph of a partial function $f:M^{|x|-1}\rightarrow M$. In particular, it follows that the case $|x|=1$ of the proof is immediate. Hence onwards we assume that $|x|>1$. 

Since, the every $b\in M^{|y|}$, the formula $\varphi(A;b)$ defines the graph of a partial function $M^{|x|-1}\rightarrow M$, for any $c\in M^{|x|-1}$ and distinct $s,t\in A_c$ it holds that\begin{equation}\label{eqn:intersec.}
\varphi(c,s; M) \cap \varphi(c,t; M) =\emptyset.
\end{equation}
Let 
\[
E=\{ a\in A : \dim \varphi(a;M)=|y|\}.
\]
By Equation~\eqref{eqn:intersec.} and the Fiber Lemma for o-minimal dimension~\cite[Chapter 4, Proposition 1.5]{dries98} we have that, for any $c\in M^{|x|-1}$, the fiber $E_c=\{ t : (c,t)\in E\}$ is finite. In particular $\dim E \leq |x|$. 

Consider a cell decomposition $\Dd$ of $M^{|x|}$ compatible with $\{A, E\}$. We apply Lemma~\ref{lem:size}, with $\Xx=\{\varphi(A;b) : b\in M^{|y|}\}$ and the $\Xx_i$ being the families $\{\varphi(D;b) : b\in M^{|y|}\}$ for $D\in \Dd$ with $D\subseteq A$, and derive that to prove the proposition it suffices to show that 
\begin{equation}\label{eqn2}
|\{\varphi(D;b) : b\in M^{|y|}\}| \leq |A|
\end{equation}
for every $D\in \Dd$ with $D\subseteq A$.
Let us fix one such cell $D$. Since otherwise the result is trivial we may assume that $D$ is infinite. Set $l := \dim D>0$.

If $l < |x|$ then there exists an injective projection $\pi_D:D\rightarrow M^{l}$. The result then follows by induction hypothesis, applied to the set $\pi_D(D)$ and formula $\varphi(\pi_D^{-1}(z); y)$ in place of $A$ and $\varphi(x;y)$ respectively. 

Now suppose that $l = |x|$. In particular, since $\dim E < |x|$, it holds that that $D\cap E = \emptyset$. 
\begin{claim}\label{claim:Ba}
For every $a \in D$ it holds that
\[
|\{\varphi(D; b) : b \in \varphi(a;M)\}| \leq |A|. 
\]
\end{claim}
\begin{claimproof}
Let $a\in D$. Since $D\cap E = \emptyset$ we have that $a\notin E$ and so
\[
\dim \varphi(a; M) < |y|.
\]
By cell decomposition we may partition $\varphi(a; M)$ into finitely many cells $\Cc$ and in order to prove the claim it suffices to show that, for every cell $C\in \Cc$, it holds that $|\{\varphi(D; b) : b \in C\}| \leq |A|$. Let us fix a cell $C\in \Cc$. Since otherwise the result is immediate we assume that it is an infinite cell. Set $d := \dim C$. Then $0 < d < |y|$.

Let $\pi_C:C \rightarrow M^d$ be an injective projection. Since $d < |y|$, we apply the induction hypothesis to the set $D$ and formula
$\varphi(x;\pi_C^{-1}(z))$
to complete the proof of the claim.
\end{claimproof}

Finally note that 
\[
\{\varphi(D;b) : b\in M^{|y|}\} \subseteq \{\emptyset\} \cup \{ \varphi(D;b) : b\in \varphi(a;M), \, a\in D\}
\]
and so by Claim~\ref{claim:Ba} we conclude that inequality~\eqref{eqn2} holds.
\end{proof}


We now derive some corollaries of Proposition~\ref{prop:main}.
Recall the preliminaries on types from Section~\ref{sec:dfns}
We use the following fact. 

\begin{fact}[O-minimal Uniform Definability of Types~\cite{cubides-ye}] \label{fact:UDT}
Let $M$ be an o-minimal structure and $\varphi(x;y)$ a formula. Then there exists another formula $\psi(y;z)$ such that, for every definable type $p(x)\in S_x^{\df}(M)$, there exists $c\in M^{|z|}$ such that 
\[
\{b \in M^{|y|} : \varphi(x;b) \in p(x)\} = \psi(M,c).
\]
\end{fact}


It is easy to see that, in any structure $M$, the space $S^{\df}(M)$ of definable types over $M$ is bounded in size by $|M|^{|L|+\aleph_0}$. For types over a definable set $A$ we may refine this bound as follows.

\begin{corollary}
Let $M$ be an o-minimal structure in a language $L$. Let $A\subseteq M$ be a definable set with $|A|>1$. Then  
\begin{equation}\label{eqn:cor1}
|S^{\df}(A)| \leq |A|^{|L|+\aleph_0}.
\end{equation}
Suppose that $A$ is infinite and let us fix a formula $\varphi(x;y)$ of $L$. For any type $p(x)\in S_x^{\df}(M)$, let $p|_{\varphi, A}(x)$ denote the restriction of $p(x)$ given by
\[
p|_{\varphi, A}(x) =  \{\varphi(x; a) : a\in A^{|y|},\, \varphi(x; a)\in p(x) \}.
\]
Then we have that 
\begin{equation}\label{eqn:cor2}
|\{ p|_{\varphi, A}(x) : p(x) \in S_x^{\df}(M)\}| \leq |A|. 
\end{equation}
\end{corollary}
\begin{proof}
Let $F$ denote the set of partitioned formulas of $L$. Note that $|F| = |L|+\aleph_0$. For any tuple of variables $x$ let $F(x)\subseteq F$ denote the set of formulas of the form $\varphi(x;y)$ for any $y$.

For any formula $\psi(y;z)$ let $\Def_\psi=\{\varphi(A^{|y|};c) : c\in M^{|z|}\}$. By Proposition~\ref{prop:main} note that $|\Def_\psi| \leq |A|+\aleph_0$.
Observe that every type $p(x)\in S_x^{\df}(A)$ is uniquely characterized by the function 
\begin{align*}
F \rightarrow \cup \{ \Def_\psi : \psi \in F\} \cup \{\emptyset\} 
\end{align*}
which maps every formula in $F\setminus F(x)$ to the empty set and, for every $\varphi(x;y)\in F(x)$, does 
\[
\varphi(x;y) \mapsto \{ a \in A^{|y|} : \varphi(x;a) \in p(x)\}.
\]

Applying Proposition~\ref{prop:main} we have that
\[
\left|\cup \{ \Def_\psi : \psi \in F\}\right| \leq |F|(|A|+\aleph_0)=|F||A|.
\]
Recall that for any infinite cardinal $\kappa$ it holds that $\kappa^\kappa = 2^\kappa$, and in particular $|F|^{|F|}=2^{|F|}$. We conclude that
\[
|S^{\df}(A)| \leq (|F||A|)^{|F|} = |F|^{|F|} |A|^{|F|} = |A|^{|F|} = |A|^{|L|+\aleph_0},
\]
proving inequality~\eqref{eqn:cor1} in the Corollary. 

Now suppose that $A$ is infinite. To prove inequality~\eqref{eqn:cor2} we apply Fact~\ref{fact:UDT} as follows. Let us fix $\varphi(x;y)$ in $F$ and let $\psi(y;z)$ be as given by Fact~\ref{fact:UDT}. Then, by Proposition~\ref{prop:main}, we have that 
\[
|\{ p|_{\varphi, A}(x) : p(x) \in S_x^{\df}(M)\}| \leq |\{\psi(A^{|y|},c) : c\in M^{|z|}\}| \leq |A|.
\] 
\end{proof}

\begin{definition}
Given a structure $M$, a topological space $(A,\tau)$ is definable in $M$ when $A\subseteq M^n$ is a definable set and there exists a definable basis for the topology $\tau$. 
\end{definition}

For a review of various classes of definable topological spaces in o-minimal structures, including definable metric spaces and manifold spaces, see~\cite{andujar_thesis}.

Recall that the weight of a topological space is the minimum cardinality of a basis for its topology. A topological space is second-countable when it has countable weight. 

The following two results follow directly from Proposition~\ref{prop:main}.

\begin{corollary}
Let $(A,\tau)$ be an infinite topological space definable in an o-minimal structure. Then $(A,\tau)$ has weight at most $|A|$. 
\end{corollary}

\begin{corollary}
The following topological spaces, which are known to be countable but not second-countable, are not definable in any o-minimal structure:
\begin{enumerate}
    \item Arens space~\cite{arens}.
    \item Maximal compact topology~\cite[Counterexample 99]{counterexamples}.
    \item One Point Compactification of the Rationals~\cite[Counterexample 35]{counterexamples}.
    \item Appert space~\cite[Counterexample 98]{counterexamples}.
\end{enumerate}
\end{corollary}


\bibliography{o-min-card}
\bibliographystyle{alpha}

\end{document}